\documentstyle[12pt]{amsart}

\input colordvi

\makeatletter
\def\LaTeX{\leavevmode L\raise.42ex
\hbox{\kern-.3em\size{\sf@size}{0pt}\selectfont A}\kern-.15em\TeX}
\makeatother

\newcommand{\BibTeX}{{\rm B\kern-.05em{\sci\kern-.025emb}\kern-.08em\TeX}}

\makeatletter
\def\@currentlabel{2.1}\label{e:dispaa}
\def\@currentlabel{2.21}\label{e:dispau}
\def\@currentlabel{2.22}\label{e:dispav}
\def\@currentlabel{2.23}\label{e:dispaw}
\def\@currentlabel{2.24}\label{e:dispax}
\def\theequation{\thesection.\@arabic\c@equation}

\newcommand{\RR}{ I\!\!R}

\newcommand{\R}{{\Bbb R}}

\newcommand{\N}{{\Bbb N}}

\newcommand{\beq}{\begin{equation}}
\newcommand{\eeq}{\end{equation}}
\newcommand{\cuad}{{\sqcap\kern-.68em\sqcup}}

\newcommand{\equ}[1]{(\ref{#1})}

\renewcommand{\theequation}{\thesection.\arabic{equation}}

\newtheorem{teo}{Theorem}[section]
\newtheorem{prop}{Proposition}[section]
\newtheorem{lemma}{Lemma}[section]
\newtheorem{corollary}{Corollary}[section]
\newtheorem{remark}{Remark}[section]
\newcommand{\bremark}{\begin{remark} \em}
\newcommand{\eremark}{\end{remark} }

\newcommand{\M}{{\cal M}}

\newcommand{\ve}{\varepsilon}

\begin{document}
\title[Eigenvalues]{  Eigenvalues for radially symmetric non-variational fully nonlinear 
operators}
%
\author{Maria J. Esteban}
\address{M.J. Esteban, Ceremade UMR CNRS 7534, Universit\'e Paris Dauphine, 75775 Paris Cedex 16, France. {\sl (esteban@ceremade.dauphine.fr)}}
\author{Patricio Felmer}
\address{\noindent P. Felmer, Departamento de
Ingenier\'{\i}a  Matem\'atica\\ and
Centro de Modelamiento Matem\'atico, UMR2071 CNRS-UChile\\
 Universidad de Chile, Casilla 170 Correo 3, Santiago, Chile. {\sl  (pfelmer@dim.uchile.cl)}}
\author{Alexander Quaas}
\address{Departamento de  Matem\'atica, Universidad Santa Mar\'{i}a
Casilla: V-110, Avda. Espa\~na 1680, Valpara\'{\i}so, Chile. {\sl (alexander.quaas@usm.cl)}
}


\begin{abstract}
In this paper we present an elementary  theory about the existence of 
eigenvalues for fully nonlinear radially symmetric 1-homogeneous operators.
 A general theory for first eigenvalues and eigenfunctions of 1-homogeneous fully nonlinear operators exists in the framework of viscosity solutions. 
Here we want to show that for the radially symmetric operators (and one dimensional) a much simpler theory can be established, and that the complete set of 
eigenvalues and eigenfuctions characterized by the number of zeroes can be
obtained.
\end{abstract}
\date{}\maketitle

\setcounter{equation}{0}
\section{Introduction}\label{sec:intro}

A fundamental step in the analysis of nonlinear equations is the understanding 
of the associated eigenvalue problem. In the case of our interest 
the question is  the 
existence of nontrivial solutions $(u,\lambda)$ of the boundary value problem
\begin{eqnarray}
F(D^2u,Du,u,x)&=&\lambda\, u \quad\mbox{in}\quad \Omega\label{E1}\\
u&=&0\quad\mbox{on}\quad \partial\Omega,\label{E2}
\end{eqnarray}
where $F$ is a positively homogeneous elliptic operator and $\Omega$ is 
a bounded smooth domain in $\R^N$, $N\ge 1$.

There is a well established theory for the first eigenvalue and eigenfunction
for this problem in the framework of viscosity solutions.
The first result in this direction is due to P.L. Lions who proved
existence of a first eigenvalue and eigenfunction for the Bellman equation in
\cite{lions1} and for the Monge-Amp\`ere equation in
\cite{lions2} by means of probabilistic arguments. 
More recently  Quaas and Sirakov addressed, by purely partial differential equations arguments, 
the general case in \cite{quaassirakov} for the existence and qualitative 
theory and the existence of solutions to the 
associated forced Dirichlet problem 
when $\lambda$ stays below the first eigenvalue. Results in this direction were
also obtained by
 Armstrong 
 \cite{Amstrong} and Ishii and Yoshimura \cite{Ishii-}.  While in \cite{quaassirakov}
convexity of $F$ is required, in \cite{Ishii-} and  \cite{Amstrong} this
hypothesis is not necessary.
Earlier partial results were obtained by
Felmer and Quaas \cite{FQ} and Quaas \cite{Q}, see also the detailed 
bibliography contained in \cite{quaassirakov}. 
Based on the eigenvalue theory just discussed, it is possible to build on 
the existence of positive (or negative) solution of the equation
\begin{eqnarray}
F(D^2u,Du,u,x)&=&\lambda\, u +f(x,u)\quad\mbox{in}\quad \Omega\label{E11}\\
u&=&0\quad\mbox{on}\quad \partial\Omega,\label{E12}
\end{eqnarray}
by means of bifurcation theory, 
using the ideas of Rabinowitz \cite{rabinowitz1},  \cite{rabinowitz2} and 
 \cite{rabinowitz}.

A better understanding of the solutions of equations \equ{E11}-\equ{E12} can
be obtained if further eigenvalues and eigenfunctions are known for 
\equ{E1}-\equ{E2}, however this has been elusive in these general 
fully nonlinear setting, except
in some particular cases in presence of radial symmetry as in the work
by Berestycki \cite{berestycki}, Arias and Campos \cite{ariascampos} for the Fucik operator,
by Busca, Esteban and Quaas \cite{BEQ} for the Pucci operator and more recently
for a more general class of extremal operators by Allendes and Quaas 
\cite{allendes-quaas}. More precisely, in \cite{BEQ}, \cite{allendes-quaas}
and \cite{ariascampos}
a sequence of eigenvalues and eigenfunctions characterized by 
their number of zeroes is constructed and a global bifurcation theory 
is obtained upon  them.

The aim of this article is to prove the existence of a sequence of 
eigenvalues and eigenfunctions for a general fully nonlinear operator 
in the radially
symmetric case, and in a self contained fashion, based on elementary arguments. 
This construction is based on the existence of two 
{\it semi eigenvalues} associated to positive and 
negative eigenfunctions in the ball and in concentric annuli, put together
via  degree theory through a Nehari type 
approach \cite{nehari}. While the spectral theory for a ball and annuli can be 
obtained as particular cases of the general results in \cite{quaassirakov},  
\cite{Amstrong} and
  \cite{Ishii-}, the general arguments to obtain the 
existence of semi eigenvalues and 
positive (negative) eigenfunctions are quite sophisticated, based on the whole
theory of viscosity solutions. When dealing with the radially symmetric
problem
in the ball or an annulus, much simpler arguments can be given. In fact, 
this is a
purely ordinary differential equation problem, which only gains some 
difficulties
in the case of a ball due to the singularity at the origin, 
but of a different order when compared with the general case.
It is our purpose to provide a simple, self contained spectral theory.

Now we present in precise terms our main theorem. On the operator $F$ we assume
the same general hypotheses as in \cite{quaassirakov}, namely:
$F:S_N\times \R^N\times\R\times B_R \to\R$, is a  continuous function, where  
$B_R$ is the ball of radius $R$, centered at the
origin and $S_N$ is the set of all symmetric $N\times N$ matrices. 
On $F$ we will make the following assumptions 
\begin{itemize}
\item[(F1)]  $F$ is  
positively homogeneous of degree $1$, that is, for all $s\ge 0$ and  
for all $(M,p,u,x)\in S_N\times \R^N\times\R\times\Omega$,
$$
F(sM,sp,su,x)=sF(M,p,u,x).
$$
\item[(F2)] There exist numbers $\Lambda\ge \lambda>0$ and $\gamma,\delta>0$ such that for all
$M,N\in S_N$, $p,q\in\R^N$, $u,v\in\R$, $x\in\Omega$
\begin{eqnarray}
\M_{\lambda,\Lambda}^-(M-N)-\gamma|p-q|-\delta|u-v|
\le F(M,p,u,x)\nonumber\\
~~~~~~~~ -F(N,q,v,x)\le 
\M_{\lambda,\Lambda}^+(M-N)+\gamma|p-q|+\delta|u-v|.\nonumber
\end{eqnarray}
Here $\M_{\lambda,\Lambda}^+$ and $\M_{\lambda,\Lambda}^-$ are the maximal and minimal Pucci operators with parameters $\lambda $ and $\Lambda$, respectively. 

\item[(F3)]   For all $M,N\in S_N$, $p,q\in \R^N$, $u,v\in\R$, $x\in\Omega$,
\begin{eqnarray}
-F(N-M,q-p,v-u,x)&\le &F(M,p,u,x)-F(N,q,v,x)\nonumber\\
&\le & F(M-N,p-q,u-v,x).\nonumber
\end{eqnarray}
\end{itemize}
The last assumption (F3) together with (F1) implies that $F$ is convex in
$(M,p,u)$, important property that we will use repeatedly in the sequel. 

In this article  we consider the extra assumption that 
the operator is radially invariant. For stating this, consider
a smooth radially symmetric function $u=u(r)$, then we have
$$
Du(x)=\frac{x}{r}\,u'(r)\quad\mbox{and}\quad
D^2u(x)=\frac{u'(r)}{r}I+\Big(u''(r)-\frac{u'(r)}{r}\Big)\frac{x\otimes x}{r^2}.
$$
Writing $m=u''(r)$ and  $p=u'(r)$, we assume
\begin{itemize}
 \item[(F4)] The operator $F$ is radially invariant, that is,
$$
F\Big(\frac{p}{r}\,I+(m-\frac{p}{r})\,\frac{x\otimes x}{r^2},\,\frac{p}{r}\,x,\,u,\,x\Big)
$$ 
 depends on $x$ only through $r$.
  \end{itemize} 

Now we can state our main theorem
\begin{teo}\label{multipleNdim}
Under assumptions (F1)-(F4), the eigenvalue problem \equ{E1}-\equ{E2} in
the ball
$B_R$
possesses  sequences of 
classical radially symmetric solutions $\{(\lambda^\pm_n, u^\pm_n)\}$,  both $u^+_n$ and $u^-_n$
with $n$ interior zeros 
$r_1<...<r_n$ and $u_n^+$ (resp. $u_n^-$) is positive (resp. negative) in 
the interval $(0, r_1)$.
Moreover the sequences $\{\lambda^\pm_n\}$ are increasing and 
the sequence $\{(\lambda^\pm_n, u^\pm_n)\}$ are complete in the sense that 
there are no  radially symmetric eigenpairs of \equ{E1}-\equ{E2} outside of them.
\end{teo}
As we already mentioned we prove this theorem relying on ordinary differential arguments. At a first step we study of the eigenvalue problem in an annulus which
becomes a regular ordinary differential equations problem. In doing so we
prove a one dimensional version of our main theorem whose precise statement is given in Theorem 
\ref{maintheorem}. The proof of the theorem uses classical existence theory
for the initial value problem together with maximum and comparison principles
obtained by means of the Alexandrov-Bakelman-Pucci (ABP) inequality. This allows to prove an existence and uniqueness theorem for a Dirichlet boundary value
problem upon which we set up a parameterized fixed point problem where Krein
Rutman theorem, in the version used by Rabinowitz \cite{rabinowitz} can be applied.
Thus we obtain a spectral theory for the first positive and negative eigenvalues in an interval, which applies also to the annulus in the radially symmetric
$N$-dimensional case.

In order to obtain the whole set of eigenvalues and eigenfunctions we use
a Nehari approach via degree theory. In this respect we notice that a qualitative property needed to use this approach is the monotonicity
of the positive and negative semi eigenvalues with respect to the interval. 
This property is a very easy consequence of a  min-max 
definition of the eigenvalues. Here we do not start in this way, but we obtain
the eigenvalues through nonlinear bifurcation theory with the help of  Krein Rutman theorem. It is interesting to see, and probably was not known, that 
the monotonicity
property can be obtained by further analyzing the proof of 
Krein Rutman theorem, see Corollary \ref{ccc}.

As a second step in the proof of Theorem \ref{multipleNdim} we study 
the eigenvalue problem in a ball, following a similar approach as in the
one dimensional case, but studying in detail the
 singularity at the origin. Regularity and compactness properties are proved
for solutions of this ordinary differential equation using simple arguments.

The paper is organized as follows. After we prove some auxiliary results in 
Section \ref{assumptions}, we treat the case of the principal eigenvalue for 
1-dimensional problems in Section \ref{principal} and we prove some qualitative properties of the eigenvalues. In Section 
\ref{eigenvalues} we prove the existence of a complete sequence of eigenvalues
and eigenfunctions in the one dimensional case. Finally, in Section 
\ref{multi} we extend the results to the radially symmetric 
multidimensional case.

\bigskip

\setcounter{equation}{0}
 \section{The unidimensional case: preliminaries}
\label{assumptions}

In this section we assume that the operator $F$ satisfies  hypotheses
(F1), (F2) and (F3) with $N=1$ and we prove a preliminary result that essentially says that we can isolate the second derivative from the equations, allowing
to use ordinary differential equations arguments to follow. We end the section with the maximum and comparison principles in this one dimensional setting.

Before continuing let us observe that, in particular, we are assuming that 
$F:\R^3\times [a,b]\to\R$ is a continuous function and it satisfies
\begin{itemize}
\item[(F2)] There are constants $\Lambda\ge\lambda>0$, $\gamma>0$ and $\delta>0$ so that
for all $(m,p,u,t),\; (m',p',u',t)\in\R^3\times [a,b]$,
\begin{eqnarray*}
-\delta|u-u'|-\gamma|p-p'|+\lambda(m-m')^+-\Lambda(m-m')^-&\le &\\
F(m,p,u,t)-F(m',p',u',t)&\le&\\
\Lambda(m-m')^+ -\lambda(m-m')^-+\gamma|p-p'|+\delta|u-u'|.
\end{eqnarray*}
\end{itemize}
Here and in what follows we write $x^+=\max\{x,0\},\; 
x^-=\max\{-x,0\}$ so that $x=x^+-x^-$.

\medskip

In the one dimensional setting the main goal of this paper is to 
study  the eigenvalue problem
\begin{equation}\label{VP} F(u'',u',u,t)=-\mu u, \quad\mbox{in}\quad [a,b],\quad u(a)=u(b)=0\end{equation}
and the auxiliary Dirichlet problem
\begin{equation}\label{Dir} F(u'',u',u,t)=f(t), \quad\mbox{in}\quad [a,b],\quad u(a)=u(b)=0,\end{equation}
In what follows we denote by $C_2(a,b)$ the space  $C^2(a,b)\cap C^1([a,b])$
and we say that $u$
is a solution of problems  \eqref{VP} and \eqref{Dir} if $u\in C_2(a,b)$ and if it satisfies the corresponding equation in $(a,b)$, together with
the boundary conditions. We notice that with our definition, a solution 
always has well defined derivatives at the extremes
of the interval $(a,b)$.

Our first result allows us to isolate $u''$ in equations  \eqref{VP} and 
\eqref{Dir}, a very convenient fact for existence and regularity analysis.
\begin{lemma}\label{despeje}
If (F1) and (F2) hold, 
there is a continuous function $G:\R^4\to \R$ so that
$$
F(m,p,u,t)=q\quad\mbox{if and only if}\quad m=G(p,u,q,t),
$$
$G$ being Lipschitz continuous in $(p,u,q)$ and monotone increasing in $q$.
\end{lemma}
\noindent
{\bf Proof.}
Using (F2), we see that
\begin{equation}\label{onto}
\lambda m^+-\Lambda m^-\le 
F(m,p,u,t)-F(0,p,u,t)\le
\Lambda m^+ -\lambda m^-,
\end{equation}
from where it follows that, for every $(p,u,t)$ fixed, $F(\cdot,p,u,t)$ is onto $\R$. Indeed, \eqref{onto} implies that $F$ is not bounded. This, together with the continuity property, proves our claim. On the other hand,
if there are $\,m, m'\,$ so that
$$
F(m,p,u,t)=F(m',p,u,t)\,,
$$
then, from (F2) again, 
$$
\lambda (m-m')^+-\Lambda (m-m')^-\le 
0\le
\Lambda (m-m')^+ -\lambda (m-m')^-
$$
from where it follows that $m=m'$.
Thus, given $(p,u,q,t)$, there is a unique $m$ so that $F(m,p,u,t)=q$, we denote by $G(p,u,q,t)$ 
such $m$. 
This function $G$ is continuous. We also prove that it is Lipschitz continuous in the 
first three variables.
Assume that
$$
q=F(m,p,u,t)\quad\mbox{and}\quad q'=F(m',p',u',t)
$$
then from (F2) we have, in case $m\ge m'$, 
$$
q-q'\ge \lambda(m-m')-\gamma|p-p'|-\delta|u-u'|,
$$
so that
$$
0\le G(p,u,q,t)-G(p',u',q',t)\le \frac{1}{\lambda}|q-q'|+\frac\gamma\lambda|p-p'|+\frac\delta\lambda|u-u'|,
$$
and if $m<m'$, then 
$$
q-q'\le -\lambda (m-m')^-+\gamma|p-p'|+\delta|u-u'|,
$$
so that
$$
0\le G(p',u',q',t)-G(p,u,q,t)\le  \frac{1}{\lambda}|q-q'|+\frac\gamma\lambda|p-p'|+\frac\delta\lambda|u-u'|.
$$
Thus, $G$ is Lipschitz continuous in $(p,u,q)$.

\smallskip
Finally, let $q\le q'$ and $m,m'$  such that
$m=G(p,u,q,t)$ and $m'=G(p,u,q',t)$. Then, $F(m,p,u,t)=q\le q'=F(m',p,u,t)$, so that from
(F2), we have
$$
-\Lambda(m-m')^-+\lambda (m-m')^+\le F(m,p,u,t)-F(m',p,u,t)=q-q'\le 0,
$$
which implies that  $m\le m'$, proving that $G(p,u,q,t)\le G(p,u,q',t)$.
 $\Box$

\bigskip
The following is a direct consequence of Lemma \ref{despeje}. 
\begin{corollary}\label{co}
Assume that $F$ satisfies (F1) and (F2) and that $u\in C_2(a,b)$ is a nontrivial solution of
$$
F(u'',u',u,t)=-\mu u, \quad\mbox{in}\quad (a,b), \quad u(a)=0
$$
then $u'(a)\not =0$.
\end{corollary}

An important ingredient in the study of fully nonlinear problems is the 
maximum and comparison principles as expressed by the ABP inequalities.
Here we present a one dimensional version:
\begin{prop}(ABP)\label{teo1}
Assume that $u\in C_2(a,b)$ is a solution of
$$
\Lambda (u'')^+-\lambda (u'')^-+\gamma |u'|\ge -f^-\quad\mbox{ in }\quad \{u>0\},
$$
with $u(a),u(b)\le 0$, then
\begin{equation}\label{max1}
\sup_{(a,b)}u^+ \le B\,\|f^-\|_{L^1(a,b)}.
\end{equation}
On the other hand, if
$u$ is a solution of
$$
\lambda (u'')^+-\Lambda (u'')^--\gamma |u'|\le f^+\quad\mbox{ in } \;\{u<0\}
$$
with $u(a),u(b)\ge 0$, then
\begin{equation}\label{max2}
\sup_{(a,b)}u^-\le B\,\|f^+\|_{L^1(a,b)}.
\end{equation}
The constant $B$ depends on $\lambda,\gamma$ and $b-a$.
\end{prop}
\noindent
The proof of this proposition can be obtained from the general 
$N$-dimensional case, see \cite{cc} for example, however in Section \S \ref{multi} we present a 
simplified proof adapted to this situation, including also the radial case. Some direct corollaries that follow from Proposition \ref{teo1} are:
\begin{corollary}
Assume that
  $F$ satisfies (F2), $\, F(m,p,u,t)$ is decreasing in $u$ and  $u\in 
C_2$. If $u$  satisfies $F(u'',u',u,t)\ge -f^-$, then \equ{max1} holds,
and if $u$  satisfies 
$F(u'',u',u,t)\le f^+$, then \equ{max2} holds.
\end{corollary}
\noindent
And the comparison principle:
\begin{corollary}\label{teo2}
Assume that $F$ satisfies (F2), (F3) and $F$ is decreasing in $u$. If $u,v\in C_2$ satisfy
$$
F(u'',u',u,t)\ge F(v'',v',v,t) \quad \mbox{ in } (a,b),\quad
$$
and $u(a)= v(a)$, $u(b)= v(b)$, then,  $u\le v$ in $[a,b]$.
\end{corollary}

\setcounter{equation}{0}
\section{A theory for the first eigenvalue and eigenfunction}
\label{principal}

The purpose of this section is to present a simplified version of the first eigenvalue theory in the unidimensional case.
We start
with an existence theorem for the Dirichlet problem in a finite interval.
\begin{teo}\label{teo3}
Assume that $F$ satisfies (F1), (F2) and (F3), then,  there exists $\kappa> 0$ such that the equation
\begin{equation}\label{dir1}
F(u'',u',u,t)-\kappa u=f(t), \quad\mbox{in}\quad (a,b),\quad u(a)=u(b)=0.
\end{equation}
has a unique solution $u\in C_2(a,b)$, for any $f\in C^0[a,b]$. 
\end{teo}
\noindent
{\bf Proof.}
First, for a given $d\in\R$, we consider the initial value problem
\begin{eqnarray*}
{ F}(u'',u',u,t)-\kappa u&=&f\;,\quad\mbox{for}\quad
t\in(a,b), \quad\\ u'(a)=d,\; u(a)=0,&&
\end{eqnarray*}
which has a unique solution since, by Lemma \ref{despeje} this equation is equivalent  to
\begin{eqnarray}
u''={ G}(u',u,f(t)+\kappa u,t) \;,&&\quad\mbox{for}\quad
t\in(a,b),\quad\label{eqq1}\\ u'(a)=d, \; u(a)=0,&&\label{eqq2}
\end{eqnarray}
with $G$ Lipschitz continuous.
We observe that the solution can be extended for all $t\in (a,b)$, since the 
nonlinearity growths less than linearly. If we denote by $u(d,t)$ the 
corresponding solution, we see that the map $\,d\mapsto u(d,b)$ is continuous.

Next, for $\,\kappa\,$ large enough, depending only on the structural 
constants of ${ F}$, 
we  consider two constants, $M_-<0<M_+\,$ such that
the constant function $\,u_+(t)=M_+\,$ satisfies
\begin{eqnarray*}
{ F}(0,0,M_+,t)-\kappa M_+&\le&f\quad\mbox{in}\quad
(a,b),\quad
\end{eqnarray*}
and the constant function $u_-(t)=M_-$ is a solution of
\begin{eqnarray*}
{ F}(0,0,M_-,t)-\kappa M_-&\ge&f\quad\mbox{in}\quad
(a,b),\quad
\end{eqnarray*}
Now we claim that there are numbers $d_1\in \R$ and $t_1\in (a,b)$ 
such that $u(d_1,t) {\geq} M_+$ for all $t\in (t_1,b]$, and similarly,  
there are numbers $d_2\in \R$ and $t_2\in (a,b)$ 
such that $u(d_2,t) {\leq} M_-$ for all $t\in (t_2,b]$.
In particular $u(d_1,b)>0$ and  $u(d_2,b)<0$.
Assuming the claim for the moment, and using the continuity of $d\to u(d,b)$
we conclude to the existence of a solution of \equ{dir1}.

Now we prove the claim. Since $G$ is Lipschitz continuous, there is a constant
$L$ such that 
$$
|G(p,u,f(t)+ku,t)|\le L(|p|+|u|+1),\quad\mbox{for all}\quad t\in[a,b],
$$
so that if for some $d_1\ge 1$ we have   $|u'(t)|\le d_1$ and $|u(t)|\le d_1$, 
then
$$
|G(u'(t),u(t),f(t)+ku(t),t)|\le 3Ld_1.
$$
Using the equation in the form \equ{eqq1}-\equ{eqq2}, we see then
 that for { $a\leq t\le
t_1:=a+ 1/4L$} we have $u'(t)\ge d_1/4$. Now we choose $d_1$ large enough so that
{ $d_1\,t_1/4>M_+$} and we find that $u(t_1)>M_+$.{ Using the Comparison Principle as given in Corollary \ref{teo2}, we see that $u(t)\ge M_+$}
for all $t\in (t_1,b]$. The case with $M_-$ is similar.
$\Box$

\medskip

Next we present an existence result that will be used in an approximation
procedure in the multidimensional radial case in Section  \ref{multi}.
\begin{teo}\label{lema1}
Assume that $F$ satisfies (F1), (F2) and (F3), then  there exists $\kappa> 0$ 
such that for every $c\in [a,b)$ and for any $f\in C^0[a,b]$ the equation
\begin{equation}\label{dir}
F(u'',u',u,t)-\kappa u=f(t), \,\,\mbox{in}\,\, (a,b),\,\, u'(c)=u(b)=0.
\end{equation}
has a unique solution $u\in C_2(a,b)$.
\end{teo}
\noindent
{\bf Proof.}
For a given $d\in\R$, we consider the initial value problem
\begin{eqnarray*}
&&u''={ G}(u',u,f(t)+\kappa u,t) \;,\quad\mbox{for}\quad
t\in(a,b),\\ && u'(c)=0, \; u(c)=d.
\end{eqnarray*}
We denote by $u(d,t)$ the 
corresponding solution and we observe that the map $\,d\mapsto u(d,b)$
is continuous.

Next, for $\,\kappa\,$ large enough (depending only on the structural 
constants of ${ F}$), 
we  consider two constants, $M_-<0<M_+\,$ such that
the constant function $\,u_+(t)=M_+\,$ satisfies
$$
{ F}(0,0,M_+,t)-\kappa M_+\le f\quad\mbox{in}\quad
(a,b),\quad u_+'(c)=0, \; u_+(b)>0,
$$
and the constant function $u_-(t)=M_-$ is a solution of
$$
{ F}(0,0,M_-,t)-\kappa M_-\ge f\quad\mbox{in}\quad
(a,b),\quad u_-'(c)=0, \;u_-(b)<0.
$$
Now we claim that for 
$d_1>M_+$ the function $u_1(t):=u(d_1,t)$ satisfies
$$
u_1(t)\ge M_+ \quad\mbox{for all }t\in {(c,b)}\,, 
$$
while for $d_2<M_-$,  $u_2(t)=u(d_2,t)$ satisfies
$$
u_2(t)\le M_- \quad\mbox{for all }t\in {(c,b)}\,.
$$
In particular $u_1(b)>0$ and  $u_2(b)<0$.
Assuming the claim for the moment, and using the continuity of $d\to u(d,b)$
we conclude to the existence of a solution of \equ{dir}.

In order to prove the claim we use the Comparison Principle as stated in Corollary \ref{teo2}. What we do is to formulate a problem in
the interval $(2c-b,b)$ by reflecting the corresponding elements. We start 
defining
$$
{ F}_c(x,y,z,t)={ F}(x,y,z,t) \quad \mbox{if}\quad t\in[c,b]
$$
and
$$
{ F}_c(x,y,z,t)={F}(x,y,z,2c-t) \quad 
\mbox{if}\quad t\in[2c-b,c].
$$
We  also reflect the solution $u(d,t)$, the right hand side $f$
 and the super and sub-solutions
$u_+$ and $u_-$. Applying the Comparison Principle contained in Corollary
 \ref{teo2} { to $\,F_c$ } we prove our claim.
Applying the same principle we also see that the solution thus found is unique,
completing the proof of the lemma.
$\Box$

\medskip

Now that we have completed the 'linear' theory we address the existence of the
first eigenvalue and eigenfunction as an application of Krein-Rutman theorem
in a form proved by Rabinowitz in \cite{rabinowitz}, see also \cite{FQ}. 
This approach will also allow us to obtain comparison results for the 
first eigenvalue
depending on the domain.

\begin{teo}\label{firstev}
Under assumptions (F1), (F2) and (F3), the eigenvalue problem 
\begin{equation}\label{eigenv}
F(u'',u',u,t)=-\mu u, \quad\mbox{in}\quad (a,b),\quad u(a)=u(b)=0
\end{equation}
has a solution  $(u^+,\lambda^+)$,  with $u^+>0$ in $(a,b)$  and another solution $(u^-,\lambda^-)$ with $u^-<0$ in $(a,b)$. Moreover, every positive 
(resp. negative) solution of equation \equ{eigenv} 
is a multiple of   $u^+$ (resp. $u^-$). 
\end{teo}
\noindent
{\bf Proof.} 
We define $K=\{u\in C[a,b]\,/\,u\ge 0, u(a)=u(b)=0\}$ and use 
 Theorem \ref{teo3} to solve Dirichlet problem 
\begin{equation}\label{dir2}
F(u'',u',u,t)-\kappa u=-g(t), \quad\mbox{in}\quad (a,b),\quad u(a)=u(b)=0.
\end{equation}
for  $g\in K$, 
provided $\kappa$ is large enough. 
We denote this solution by  ${\cal L}(g)$ and define the operator
$T: \R^+\times K\to K$ as $T(\mu,f)=\mu {\cal L}(f)$. The operator
$T$ is well defined and, {as a consequence of Corollaries \ref{co} and 
\ref{teo2}}, {$T(\mu,f)>0$} for
every {$f\in K\setminus\{0\}$}, $\mu>0$. Moreover $T$ is compact and $T(0,g)=0$ for
every $g\in K$, so it satisfies the hypothesis to obtain the existence of a family  fix points see 
Corollary 1 of Theorem VIII 1 in \cite{rabinowitz}. Notice that  this Corollary 1 is the main argument in proving Krein-Rutman Theorem. 

Take $u_0\in K\setminus\{0\}$, then 
there exists $M>0$ such that $M {\cal L} (u_0)\ge u_0$, since the 
contrary is not compatible with Corollary \ref{co}.
Define now
$T_\varepsilon : \RR^+\times K \to K$
as
$T_\ve(\mu,u)= \mu{\cal L}(u)+\mu\varepsilon {\cal L}(u_0),$ for $\ve>0$.
Then, from Corollary 1 in \cite{rabinowitz}, there exists an unbounded
 connected 
component ${\cal C}_\varepsilon$
 of
solutions to $T_\varepsilon (\mu,u)=u$,
moreover  ${\cal C}_\varepsilon\subset [0,M]\times K$.
To see this fact, let $(\mu,u)\in {\cal C}_\varepsilon$, then
\begin{eqnarray*}
u=\mu {\cal L} (u)+\mu\varepsilon {\cal L}  (u_0)\,.
\end{eqnarray*}
Hence
$u\ge\mu\varepsilon {\cal L}(u_0)\ge\frac{\mu}{M}\varepsilon u_0.
$
If we apply ${\cal L}$ we get
\begin{eqnarray*}
{\cal L} (u)\ge\frac{\mu}{M}\varepsilon {\cal L} (u_0)\ge\frac{\mu}{M^2}
\varepsilon u_0.
\end{eqnarray*}
But $u\ge\mu {\cal L} (u)$, then $u\ge (\frac{\mu}{M})^2
\varepsilon u_0$. By recurrence we get
\begin{eqnarray*}
u\ge(\frac{\mu}{M})^n\varepsilon u_0\quad\mbox{for all}\quad n\geq 2
\end{eqnarray*}
and we conclude that $\mu\le M$. This and the fact that 
 ${\cal C}_\varepsilon$ is unbounded implies that  there exists
$(\mu_\ve,u_\ve)\in {\cal C}_\varepsilon$
such that $\|u_\ve\|_\infty=1$.{ This and Theorem \ref{lema1} imply a uniform bound in $C_2(a,b)$, allowing us to  pass to the limit as $\ve\to 0$ to find}
$\mu^+ \in[0,M]$ and $ u^+>0$
such that $u^+=\mu^+  {\cal L} (u^+)$. From here we also deduce that
$\mu^+>0$ and then we define $\lambda^+=-\kappa+\mu^+$. 
For the simplicity and the isolation of the eigenvalue, we can use an argument similar to the
one given in  \cite{rabinowitz}, {later} adapted to a situation like ours in 
 \cite{FQ}.
The same argument can be applied to $-F(-m,-p,-u,t),$  a concave operator, to
obtain $(u^-,\lambda^-)$. $\Box$

\medskip

In what follows we denote by $\lambda^+(t_1,t_2)$ the first eigenvalue associated
to a positive eigenfunction, and $\lambda^-(t_1,t_2)$  the first 
eigenvalue associated
to a negative eigenfunction, 
given 
in 
Theorem \ref{firstev} for the problem  \equ{eigenv} in the interval $(t_1,t_2)\subset (a,b)$. 
\begin{corollary}\label{ccc}
If $(a_1,b_1)\subset (a,b)$ and $(a_1,b_1)\not = (a,b)$
then 
$$
\lambda^{\pm}(a_1,b_1)> \lambda^\pm(a,b).
$$
\end{corollary}
\noindent
{\bf Proof.}
We consider the eigenpair $(\lambda_1^+,u_1^+)$, $\mu_1^+=\kappa+\lambda_1^+$,
 given by Theorem \ref{firstev}
on the interval $(a_1,b_1)\subset (a,b)$, so that
$$
F((u_1^+)'', (u_1^+)',u_1^+,t)=-\lambda_1^+u_1^+\quad\mbox{in}\quad (a_1,b_1).
$$  
If $\underbar u$ is the function obtained by extending $u_1^+$ by zero  to 
 the whole interval $[a,b]$, we define $\tilde u$ as the unique solution of
$$
F(\tilde u'', \tilde u',\tilde u,t)-\kappa \tilde u=-\mu_1^+\underbar u \quad\mbox{in}\quad (a,b),\quad u(a)=u(b)=0.
$$
Then, using Comparison and Strong Maximum Principle in the interval 
$[a_1,b_1]$ we see that $\tilde u>u_1^+$ in $[a_1,b_1]$ and consequently $\tilde
 u>\underbar u$ in $(a,b)$.
If we define $w={\cal L}(\underbar u)$ and $v={\cal L}(\tilde u)$, then
$$
F(w'',w',w,t)-\kappa w=-\underbar u>-\tilde u=F(v'',v',v,t)-\kappa v,
$$
so, again by Comparison and Strong Maximum Principle, $w< v$, which implies 
$\tilde 
u=\mu_1^+{\cal L}(\underbar u)< \mu_1^+{\cal L}(\tilde u)$. Here we may replace
$\mu_1^+$ by a slightly smaller value, without changing the strict inequality. 
Now, 
the application of the construction
in the proof of {Theorem \ref{firstev}, with $\,u^0=\tilde u\,$ and} $\,M<\mu_1^+$, 
implies that $\mu^+:=\lambda^+(a,b)+\kappa< \mu_1^+$, thus completing the proof for $\lambda^+$. The proof
for $\lambda^-$ is similar.
$\Box$ 
\begin{corollary}\label{ultimo}
The functions $\lambda^+,\lambda^-:\{(t_1,t_2)\,\,/\,\, a\le t_1<t_2\le b\} \to \R$ are continuous and
$$
\lim_{t_2-t_1\to 0^+}\lambda^+(t_1,t_2)=
\lim_{t_2-t_1\to 0^+}\lambda^-(t_1,t_2)=\infty.
$$
\end{corollary}
\noindent
{\bf Proof.}
The continuity of these functions is a consequence of the 
uniqueness of the eigenvalues for positive (negative) eigenfunctions. 
While the limit is a consequence of Proposition \ref{teo1}, in fact
denoting $\mu^+=\lambda+\kappa$, with $\kappa$ as in Theorem \ref{teo3},
from \equ{max1} we obtain that
$$
\sup_{(t_1,t_2)}u^+ \le B\mu^+\,\|u^+\|_{L^1(t_1,t_2)}\le  B\mu^+(t_2-t_1)\sup_{(t_1,t_2)}u^+,
$$
which completes the proof.
$\Box$

\section{Multiple eigenvalues and eigenfunctions in the unidimensional case.}\label{eigenvalues}

In this section we consider the existence of higher eigenvalues, associated to
changing-sign eigenfunctions in the general setting already defined
in Section \S 2. More precisely,  we prove the following theorem
\begin{teo}\label{maintheorem}
Under assumptions (F1), (F2) and (F3), the eigenvalue problem
\begin{equation}\label{eigenvpm}
F(u'',u',u,t)=-\mu u, \quad\mbox{in}\quad (a,b),\quad u(a)=u(b)=0
\end{equation}
has two sequences of 
solutions $\{(\lambda^\pm_n, u^\pm_n)\}$
such that $u^\pm_n$ have both $n$ interior zeros 
$t_1<...<t_n$ and $u_n^+$ (resp. $u_n^-$) is positive (resp. negative) in 
the interval $(a, t_1)$, negative (resp. positive) on $(t_1, t_2)$.
Moreover the sequence $\{\lambda^\pm_n\}$ is increasing and 
the sequence $\{(\lambda^\pm_n, u^\pm_n)\}$ is complete in the sense that 
there are no  eigenpairs of \equ{eigenvpm} {outside these sequences}. 
\end{teo}
We devote this section  to the proof of this theorem using degree theory.
We start with a given
$n\in\N$, $n\ge 1$ and  we define 
$$
\Delta_n=\{(t_1,...,t_n)\,\,/\,\, a< t_1<t_2<...<t_n< b\},
$$
$t_0=a$ and $t_{n+1}=b$ and the  function $V: \Delta_n\to\R^n$ as
$$
V_i(\vec t)=\lambda^{(-1)^i}(t_{i-1},t_i)- \lambda^{(-1)^{i+1}}(t_{i},t_{i+1}),\qquad i=1,...,n\,,
$$
where  by $\lambda^{(\pm 1)}$ we mean $\lambda^\pm$. We observe that under  
our assumptions, Corollary \ref{ultimo} implies that the 
function $V$ is continuous in $\Delta_n$. We have the following
\begin{teo}\label{functionV}
Under assumptions (F1), (F2) and (F3),
for every $n\in\N$ there is $\vec t\in \Delta_n$ such that 
\begin{equation}
V(\vec t)=0.\label{zero}
\end{equation}
\end{teo}
\noindent
{\bf Proof.} The simplex $\Delta_n$ has a boundary composed by $n+1$ faces 
$$
F_i=\{(t_1,...,t_n)\,\,/\,\, 0=t_0\le t_1\le t_2\le ...t_i=t_{i+1}...\le t_n\le t_{n+1}= 1\},
$$
for $i=0,1,...n$. Each face has an exterior normal vector $T_i\in \R^n$ given by
$T_i=e_i-e_{i+1}$ for $i=1,...,n-1$, where $e_1, e_2,..., e_n$ is the canonical
basis of $\RR^n$. The extreme cases are 
$T_0=-e_1$ and 
$T_{n}=e_n$. 

Assume now that we have a sequence $\{\vec t_k\}$ of points in
$\Delta_n$ approaching the interior 
of a face. More precisely, if  ${\rm int}(F_i)$ represents
the interior of $F_i$ relative to $F_i$, we assume that 
$\vec t_k\to {\rm int}(F_i)$ as $k\to \infty$, for some $i=0, 1, ..., n$.
 By definition of $F_i$, as $k$ goes to $+\infty$,  $(\vec t_k)_{i+1}-(\vec t_k)_{i}\to 0$, 
and then,  according to
the definition of $V$ and Corollary \ref{ultimo}, we have that 
$ V_i(\vec t_k)\to-\infty$ and $ V_{i+1}(\vec t_k)\to +\infty$ as $k$ goes to $+\infty$, while the other components of $V$ remain bounded. 

From this  discussion we conclude that for all $i=0,1,...,n+1$
$$
\lim_{\vec t\to {\rm int}(F_i)} V(\vec t)\cdot T_i=-\infty.
$$
A well-known corollary of degree theory 
proves then the existence of a zero for $V$, {showing the existence of $\,(\lambda_n^-, u_n^-)$.}

{If we observe the definition of $V$ we see that the first component $V_1$ 
is associated to $\lambda^-(t_0,t_1)$ and $\lambda^+(t_1,t_2)$, so that the 
eigenfunction that we can construct out of solutions of equation \equ{zero}
will start being negative. 
For eigenfunctions starting with
positive values in the first interval $(t_0,t_1)$ we need to define the above arguments to the slightly
modified function
$$
\tilde V_i(\vec t)=\lambda^{(-1)^{i+1}}(t_{i-1},t_i)- 
\lambda^{(-1)^{i}}(t_{i},t_{i+1}),\qquad i=1,...,n\,.
$$
$\Box$}

\medskip

\noindent
{\bf Proof of Theorem \ref{maintheorem}.}
Given a solution
$\vec t\in\Delta_n $ of \equ{zero} we proceed to construct an eigenfunction
as follows. 
On the interval $(a, t_1)$ we define $\,u^-_n\,$ as {$\,u^-(a,t_1)$}. 
Then, on $(t_1, t_2)$ the function $\,u^-_n\,$ will be equal to 
$\alpha_1u^+(t_1, t_2)$, where $\alpha_1$ is chosen so that 
$(u^-)'(a,t_1)(t_1)= \alpha_1\,(u^+)'(t_1, t_2)(t_1)$. The existence of 
$\alpha_1$ is a consequence of Corollary \ref{co}. Here we denote by 
$u^\pm(t,s)$ the corresponding positive or negative eigenfunction on 
the interval
$(t,s)$. Repeating this argument we will finally arrive to the  function 
$u_n^-$, which is of class 
$C^1[a,b]$ and of class $C^2$ in the interior of every interval of the
form $(t_i,t_{i+1})$. Then we use the equation satisfied by
each partial eigenfunction and the continuity of $F$, rather than that of $G$,
to find that $u_n^-$ is of class $C^2(a,b)$. The associated eigenvalue is simply
$\lambda_n^-=\lambda^-(a,t_1)$.

For proving uniqueness we assume we have a second eigenpair $(\lambda, v)$
associated
 with $n$ such that there exist values $a<s_1<s_2<...<s_n<b$ and $v$ changes 
sign at those points, starting with negative values in the interval $(a,s_1)$.
If $\lambda=\lambda_n^-$ then by Corollary \ref{ccc} we necessarily have
$s_i=t_i$ for all $i=1,2,...n$ and then the simplicity and isolation of the
first eigenfunctions proved in Theorem \ref{firstev} completes the argument.

Now we assume that $\lambda>\lambda_n^-$, then by  Corollary \ref{ccc} we have
$s_1<t_1$ and then 
\begin{equation}\label{desi}
\lambda>\lambda^-(0,t_1).
\end{equation} 
We either have 
$1\le i\le n-1$ such that $(t_i,t_{i+1}) 
\subset (s_i,s_{i+1})$ or $s_n\le t_n$. In the first case,
if $i$ is odd $\lambda^+(t_i,t_{i+1})\ge \lambda$ and if
$i$ is even 
$\lambda^-(t_i,t_{i+1})\ge \lambda$, contradicting \equ{desi} in both
cases. In the second case, {$\lambda\le \lambda^+(t_n,b)$}, if
$n$ is odd, contradicting \equ{desi} again and similarly if $i$ is even.
$\Box$

\setcounter{equation}{0}
 \section{ The eigenvalue and eigenfunction theory in the radial case}
\label{multi}
We devote this section to prove our main theorem. We assume that $N>1$ and that
the operator $F$ satisfies (F1), (F2), (F3) and it is radially invariant, that
is, it satisfies also (F4). Our purpose is to study the eigenvalue problem 
\equ{E1}-\equ{E2} where $\Omega=B_R$,  is the ball of radius $R$ centered at the origin.

We start with some notation. Given our operator $F$ we define
 ${\cal F}:\R^4\times \R_+\to\R$ as
$$
{\cal F}(m,\ell,p,u,r)=F(\ell I+(m-\ell)e_1\otimes e_1,{p}e_1,u,re_1)
$$
and consider the operators 
$$
{P}^+(a,b)=\Lambda (a^++(N-1)b^+)-\lambda(a^-+(N-1)b^-)
$$
and 
$$
{P}^-(a,b)=\lambda (a^++(N-1)b^+)-\Lambda(a^-+(N-1)b^-).
$$
Here  $m$ stands for $u''(r)$, $p$ for $u'(r)$ and $\ell$ for $\frac{u'(r)}{r}$.
Under assumption (F4), we may 
rewrite  hypothesis (F2) in this radially symmetric setting as follows
\begin{itemize}
\item[(F2)] There exist $\gamma,\delta>0$ such that for all
$m,m',\ell,\ell',p,p',u,u\in\R$, $r\in[0,R]$,
\begin{eqnarray*}\nonumber
&&P^-(m-m',\ell-\ell')-\gamma|p-p'|-\delta|u-u'|\le {\cal F}(m,\ell,p,u,r)\\
&& -{\cal F}(m',\ell',p',u',r)\le
P^+(m-m',\ell-\ell')+\gamma|p-p'|+\delta|u-u'|\label{hypH2}
\end{eqnarray*}
\end{itemize}

%
%

The proof of Theorem \ref{multipleNdim} follows the general lines of that of 
Theorem \ref{maintheorem}. The new difficulty is the singularity present 
now at $r=0$. We deal with it using an approximation procedure:  
in the interval $[\ve, R)$ we apply the results of the previous sections. 
Then we obtain uniform estimates on the approximated solutions and 
their derivatives in order to pass to the limit. In the rest 
of this section we state and prove all the ingredients necessary to do this, 
and so to complete the proof of Theorem \ref{multipleNdim}. We also prove 
the (ABP) inequality for the multidimensional radial case, 
{and thus also that of Proposition \ref{teo1}}.

The next lemma is the analogue of  Lemma \ref{despeje} and it can be proved following the same arguments.
\begin{lemma}\label{DES} If (F1), (F2), (F3) and (F4) hold true, then,

\medskip
1. There is a continuous function ${\cal G}:\R^4\times \R_+\to \R$ so that
$$
{\cal F}(m,\ell,p,u,r)=q\quad\mbox{if and only if}\quad m={\cal G}(\ell,p,u,q,r)
$$
and ${\cal G}$ is Lipschitz continuous in $(\ell,p,u,q)$.

\medskip
2. There is a continuous function ${\cal G}_1:\R^3\times \R_+\to \R$
such that
$$
{\cal F}(\ell,\ell,p,u,r)=q\quad\mbox{if and only if}\quad \ell={\cal G}_1(p,u,q,r)
$$
and ${\cal G}_1$ is Lipschitz continuous in $(p,u,q)$.
\end{lemma}

\medskip

The following is a regularity result, extending the second derivative of a
solution to the origin, the only point in the domain that makes a 
difference with the one 
dimensional case.
\begin{lemma}\label{reg} Assume that (F1), (F2), (F3)) and (F4) hold true and
assume also that $f$ is
a continuous function in $[0,R]$ and $u:[0,R]\to\R$ is a solution of 
\begin{equation}\label{ed} 
{\cal F}(u'',\frac{u'}{r},u',u,r)=f(r) \quad \mbox{in}\quad (0,R)
\end{equation}
with boundary conditions
\begin{equation}\label{bc}
u'(0)=0,\; u(R)=0\,.
\end{equation}
If 
the functions $\,|u''(r)|\,$ and $\, |\frac{u'(r)}{r}|\,$ are bounded in $\,(0,R)$
then:

\medskip
1. The limit 
$$
\lim_{r\to 0} \frac{u'(r)}{r}
$$
exists
and consequently $u''(0)$ is well defined.

\medskip
2. The function $u(x)=u(|x|)$ is a $C^2(B_R)$-solution 
to the partial differential equation
\begin{equation}
F(D^2u,D u,u,x)=f  \quad \mbox{in}\quad B_R,
\end{equation}
with boundary condition
\begin{equation}
u=0  \quad \mbox{on}\quad \partial B_R.
\end{equation}
\end{lemma}
\noindent
{\bf Proof.}
 We use Lemma \ref{DES}  
to write
$$
u''={\cal G}(\frac{u'}{r},u',u,f,r)\,,
$$
and then, using the boundary condition and writing $\ell=\frac{u'}{r}$, we find
$$
r\ell=\int_0^r {\cal G}(\ell,u',u,f,s)ds\,.
$$
Differentiating the above functional equality, we get
$$
r\ell'+\ell= {\cal G}(\ell,u',u,f,r).
$$
Assume, by contradiction, that $\ell$ does not converge as $r\to 0^+$. 
Then there are two numbers $a<b$ and two sequences $\{r_n^+\}$,  $\{r_n^-\}$ such that
$$
\lim_{n\to\infty}r_n^+=\lim_{n\to\infty}r_n^-=0,
$$
and 
$$
\ell'(r_n^+)=\ell'(r_n^-)=0,\quad \lim_{n\to\infty}\ell(r_n^+)=b,\quad
 \lim_{n\to\infty}\ell(r_n^-)=a.
$$
Then we have
$$
\ell(r_n^\pm)={\cal G}(\ell(r_n\pm),u'(r_n^\pm), u(r_n^\pm ),f(r_n^\pm ),r_n^\pm )
$$
and also,
$$
\ell(r_n^\pm)={\cal G}_1(u'(r_n^\pm), u(r_n^\pm ),f(r_n^\pm ),r_n^\pm ).
$$
Since $f, u$ and $u'$ are continuous at $r=0$ as well as ${\cal G}_1$ we have that
$$
\lim_{n\to\infty }\ell(r_n^+) = \lim_{n\to\infty } \ell(r_n^-),
$$
which is a contradiction.
$\Box$

\bigskip

Next we prove (ABP) inequality in the multidimensional 
radial case. 
\begin{prop}\label{teo1radial}
Assume that $u\in C_2(0,R)$ is a solution of
$$
P^+(u'',\frac{u'}{r})+\gamma |u'|\ge -f^-\quad\mbox{ in }\quad \{u>0\},
$$
with $u(R)\le 0$, then
\begin{equation}\label{max1radial}
\sup_{(0,R)}u^+ \le B\,\|f^-\|_{L^N(B_R)}.
\end{equation}
On the other hand, if
$u$ is a solution of
$$
P^-(u'',\frac{u'}{r})  -\gamma |u'|\le f^+\quad\mbox{ in } \;\{u<0\}
$$
with $u(R)\ge 0$, then
\begin{equation}\label{max2radial}
\sup_{(0,R)}u^-\le B\,\|f^+\|_{L^N(B_R)}.
\end{equation}
The constant $B$ depends on $N, \lambda,\gamma$ and $R$.
\end{prop}
\noindent
{\bf Proof.} Let $\displaystyle{l_0=\frac{\sup_{(0,R)}u}{R}\,}$ and denote 
by $r_0$ a maximum point of $u$ in $(0,R)$. There exists a point 
$r_-\in (0,R)$ such that  $-u'(r_-)=l_0$ and
$0\leq - u'(r) \leq l_0$ in the interval $(r_0,r_-)$. Moreover, we can find 
a set $I$ (union of intervals) in $(r_0, r_-)$ so that $u''\leq 0$ in $I$ and 
$u'(I)=(0,l_0)$. We observe that on $I$ both $u''$ and $u'$ are 
 non-positive and then 
$$
P^+(u'',\frac{u'}{r})= \lambda(u''+(N-1)\frac{u'}{r})\quad\mbox{for all} \quad r\in I.
$$
Then, for any $k > 0$,
\begin{eqnarray*}
\ln(1+\frac{l_0^N}{k})&=&
\int_{0}^{l_0^N}\frac{dz}{z+k}\\
&\leq&\int_{I} \frac{-N(-u'(r))^{N-1}u''(r)\,dr}{(-u'(r))^N+k}\\
&{=}&{N}\int_{I} \left(\frac{-u'(r)}{r}\right)^{N-1}(-u''(r))\,
\frac{r^{N-1}dr}{(-u'(r))^N+k}\\
&\leq&{N}\int_{I} \left(-u''(r)-(N-1)\frac{u'(r)}{r}\right)^N
\,\frac{r^{N-1}dr}{(-u'(r))^N+k}\\
&\leq&\frac{2^NN}{\lambda^N}
 \int_{I}\left(\frac{|f^-|^N}{k}+\gamma^N\right)r^{N-1}dr\\ 
&{\le}&\frac{2^NN}{\lambda^N}\left(\frac{1}{k}\|f^-\|^N_{L^N(I)}+{\frac{(\gamma R)^N}N}
\right). 
\end{eqnarray*}
This implies that  $\|f^-\|_{{L^N(I)}}>0$, since $k$ is arbitrary. Now we choose
$k$ so that
$k\lambda^N=\|f^-\|_{L^N(I)}^N$, and find
$l_0\leq C \|f^-\|_{L^N(I)}\leq C \|f^-\|_{L^N(B_R)}$ for some constant $C>0$
 depending on $N, \lambda,\gamma$ and $R$. 
 $\Box$

\begin{remark}
When $N=1$ the proof just presented reduces to a proof of Proposition 
\ref{teo1} with
the obvious change in the domain in order  to consider a general 
interval $(a,b)$.
\end{remark}

Next corollaries follow from Proposition \ref{teo1radial}.
\begin{corollary}
Assume that hypotheses (F1)-(F4) hold and additionally that
  $\, {\cal F}(m,\ell,p,u,r)$ is decreasing in $u$. {If $u\in
C_2(0,R)$}  satisfies ${\cal F}(u'',u'/r,u',u,r)\ge -f^-$, then 
\equ{max1radial} holds,
and if {$u\in
C_2(0,R)$}  satisfies
${\cal F}(u'',u'/r,u',u,r) \le f^+$, then \equ{max2radial} holds.
\end{corollary}
\noindent
The following comparison principle also follows from Proposition \ref{teo1radial}:
\begin{corollary}\label{teo2radial}
Assume that hypotheses (F1)-(F4) hold  and additionally that
  $\, {\cal F}(m,\ell,p,u,r)$ is decreasing in $u$.
If $u,v\in C_2(0,R)$ satisfy
$$
F(u'',u'/r,u',u,t)\ge F(v'',v'/r,v',v,t) \quad \mbox{ in } (0,R),\quad
$$
and $u(R)= v(R)$, $u'(0)=v'(0)=0$, then,  $u\le v$ in $[0,R]$.
\end{corollary}

{As in Section \ref{principal}, before proving the existence of eigenvalues and eigenfunctions, we prove the existence of solutions for a related Dirichlet problem, as follows.}
\begin{teo}\label{theo5} Assume that (F1)-(F4) hold true.
There is $\kappa>0$  so that the equation
\begin{eqnarray}\label{edk}
{\cal F}(u'',\frac{u'}{r},u',u,r)-\kappa u&=&f\quad\mbox{in}\quad
(0,R),\quad\\ u'(0)=0, \; u(R)=0\,,&&\label{bck}
\end{eqnarray}
possesses a unique solution for any given continuous function $f$. 
\end{teo}
The proof of this theorem can be done  through an
approximation procedure and using only elementary ODE arguments. 
In this direction
we have the following two results.
\begin{lemma}\label{aprox} Assume assumptions (F1)-(F4).
There  is $\,\kappa>0$ (independent of $\ve$) so that for any given {$f\in C^0[0,R]$} and $\ve>0$, there exists a unique 
solution $u_\ve$ of
\begin{eqnarray}\label{edke}
{\cal F}(u'',\frac{u'}{r},u',u,r)-\kappa u&=&f\quad\mbox{in}\quad
(\ve,R),\quad\\ u'(\ve)=0, \; u(R)=0\,.&&\label{bcke}
\end{eqnarray}
\end{lemma}
\noindent
The proof of this proposition is completely similar to that of 
 Theorem \ref{lema1} so we omit it.
The following lemma provides  estimates for the solution $u_\ve$, independent of $\ve$  and its proof is inspired of that of  Lemma 2.2 in \cite{FQ}.
\begin{lemma}\label{cotas} Assume that (F1)-(F4)  hold true and
let $\,u_\ve\,$ be the solution to \equ{edke}-\equ{bcke} given by
 Lemma \ref{aprox}. Then there is a constant $C$, {independent of $\ve$,} such that
\begin{eqnarray*}
|\frac{u'_\ve(r)}{r}| \le C\quad\mbox{and}\quad |u''_\ve(r)|<C,\quad
\mbox{for all }
\ve>0, \;r\in [\ve,R].
\end{eqnarray*}
\end{lemma}
\noindent
{\bf Proof}. 
We first claim that if $u_\ve(r)$ and $u_\ve'(r)$ are uniformly
bounded in $[\ve,R]$,  then
${u'_\ve(r)}/{r}$ and $u_\ve''(r)$ are uniformly bounded in $[\ve,R]$.
By contradiction, suppose the existence of two sequences $\ve_n\to 0$ and $r_n\in (\ve_n,R]$ such that
\begin{eqnarray*}
\lim_{ n\to +\infty}\frac{u'_n(r_n)}{r_n}=-\infty,
\end{eqnarray*}
where we write $u_n=u_{\ve_n}$.
From \equ{edke}, (F2) and our assumption on $u_\ve(r)$ and $u_\ve'(r)$,
we have that
$u''_n(r_n)\to+\infty$ as $n\to+\infty$.

If $u''_n(r)>0$
for all $r\in (\ve_n,r_n]$, then $u'_n(r_n)>0$, which is impossible.
Thus,  for all $n$ there exists  $\bar r_n\in (\ve_n, r_n)$ such 
that $u''(\bar r_n)=0$ and
$u''_n(r)>0$ for all $r\in (\bar r_n,r_n)$.
Hence $u'(\bar r_n)< u'(r_n)$,
which implies that
\begin{eqnarray*}
\lim_{n\to +\infty}\frac{u'(\bar r_n)}{\bar r_n}=-\infty \quad  \mbox{ and }\quad
u''(\bar r_n)=0,
\end{eqnarray*}
which is again impossible by \equ{edke}, (F2) and our assumption on 
$u_\ve(r)$ and $u_\ve'(r)$.
Suppose next that
for a sequence of points $r_n\in (\ve_n,R)$
we have
\begin{eqnarray*}
\lim_{n\to +\infty}\frac{u'_n(r_n)}{r_n}=+\infty,
\end{eqnarray*}
then with a similar argument we also get a contradiction.
Thus, we have have that $\{{u'_\ve(r)}/{r}\}$ is bounded and as before we
conclude that 
 $\{u_\ve''(r)\}$ is bounded, proving the claim.

Suppose now that $\{\beta_\ve \}$ is unbounded with 
$\beta_\ve=\|u_\ve\|_\infty+\|u_\ve'\|_\infty$.
Define $v_\ve(r)={u_\ve(r)}/{\beta_\ve}$. 
Then $\{v_\ve\}$ and $\{v_\ve'\}$ are bounded
and  $v_\ve$ satisfies 
\begin{eqnarray}\label{edkee}
{\cal F}(v_\ve'',\frac{v_\ve'}{r},v_\ve',v_\ve,r)-\kappa v_\ve
&=&\frac{f}{\beta_\ve}\quad\mbox{in}\quad
(\ve,R),\quad\\ v_\ve'(\ve)=0, \; v_\ve(R)=0.&&\label{bckee}
\end{eqnarray}
Using the claim again, we conclude that for a positive constant $C$
\begin{eqnarray*}
|\frac{v'_\ve (r)}{r}|<C,\quad |v''_\ve (r)|<C,\quad\mbox{for all}\quad r\in 
[\ve,R],
\end{eqnarray*}
To proceed with the proof now we use the Arzela-Ascoli theorem, 
and find a sequence $v_{\ve_n}\to v$ 
uniformly in $C^1([0,R])$ to a solution $v\in C_2(0,R)$ 
of \equ{edk}-\equ{bck} with  
right-hand side equal to $0$. 
At this point we may use (ABP)
 inequality as given in Proposition \ref{teo1radial}, to obtain that this 
equation  has a unique solution by the comparison principle 
given in Corollary \ref{teo2radial},
so $v\equiv0$. But this is impossible since
$\|v_\ve\|_\infty+\|v_\ve'\|_\infty=1$ for all $\ve$.
 $\Box$

\begin{remark}
If we use the claim given in the first part of the proof of Lemma  
\ref{cotas} we see that for 
the function
$v$ defined as the limit of $v_\ve$, there is  a constant $C$ so that
\begin{eqnarray*}
|\frac{v'(r)}{r}| \le C\quad\mbox{and}\quad |v''(r)|<C,\quad
\mbox{for all }
r\in (0,R].
\end{eqnarray*}
Then we may apply Lemma \ref{reg} to find that $v''(0)$ is well defined and
thus
$v$ is a solution to the corresponding partial  differential equation 
in the ball with zero right hand side.  
\end{remark}

\noindent
{\bf Proof of Theorem \ref{theo5}}
Using Proposition \ref{theo5} we obtain a sequence of approximating solutions for
\equ{edk}-\equ{bck}. Then we use Lemma \ref{cotas} to obtain estimates that allows us
to use the Arzela-Ascoli theorem as at the end of the proof of Lemma \ref{cotas}
to obtain a solution of the problem. $\Box$

Finally we  state  a  compactness lemma, whose proof is similar to that of Lemma
\ref{reg} and which is necessary to use Krein-Rutman theory to find the
first eigenvalues. 
\begin{lemma}\label{cotass} If (F1)-(F4) hold true,
let  $u_n$ be the solution of equation \equ  {edk}-\equ{bck} with
right hand side $f_n$, where $\{f_n\}$ is a uniformly bounded sequence of continuous functions in the interval $[0,R]$. Then, there is a constant $C$, independent of $n$, such that
\begin{eqnarray*}
|\frac{u'_n(r)}{r}| \le C\quad\mbox{and}\quad |u''_n(r)|<C,\quad
\mbox{for all }
r\in (0,R].
\end{eqnarray*}
\end{lemma}
Next we have the existence of the first eigenvalues in the ball. This theorem is a particular case of the general eigenvalue theory for fully nonlinear
equations. Here we have provided a proof which relies on elementary
arguments. 
\begin{teo}\label{firstev1}
Under assumptions (F1), (F2), (F3) and (F4), the radially symmetric 
eigenvalue problem \equ{E1}-\equ{E2} in $\Omega=B_R$
has a solution  $(\lambda^+,u^+)$,  with $u^+>0$ and radially symmetric in
$B_R$ and another solution $(\lambda^-,u^-)$ with $u^-<0$ and 
radially symmetric in
$B_R$. Moreover
\begin{itemize}
\item[i)] $\lambda^+\le \lambda^-$.
\item[ii)]  Every positive
(resp. negative) solution of equation \equ{eigenv}
is a multiple of   $u^+$ (resp. $u^-$).
\item[iii)] If $\lambda^\pm(R)$ denotes the eigenvalue in $B_R$ then $\lambda^\pm(R)<\lambda^\pm(R')$ if $R>R'$.
\item[iv)] $\lambda^\pm(R)\to \infty$ if $R\to 0$.
\end{itemize}
\end{teo}
\noindent
{\bf Proof.} With the aid of Theorem \ref{theo5} and Lemma \ref{cotass} 
we can follow step by step the arguments given in the proof of Theorem 
\ref{firstev} to obtain the existence of the eigenvalues and eigenfunctions. The qualitative properties
are proved similarly as in the one dimensional case shown in Section \S3.
$\Box$

\medskip

\noindent
{\bf Proof of Theorem \ref{multipleNdim}.} 
The arguments are the same as those given in the proof of  Theorem \ref{maintheorem}. $\Box$

\bigskip

{\bf Acknowledgements:}
P.F. was  partially supported by Fondecyt Grant \# 1070314,
FONDAP and BASAL-CMM projects
 and  Ecos-Conicyt project C05E09.

A. Q. was partially supported by Fondecyt Grant \# 1070264 and USM Grant \#  12.09.17.
and Programa Basal, CMM. U. de Chile. .

\end{document}